\documentclass[]{amsart}
\usepackage{amsmath,amsthm,amssymb}
\usepackage[numbers]{natbib}
\usepackage{subfigure}
\usepackage[dvips]{graphics}
\newcommand{\hm}{\hat{\mu}^{(s)}}
\newcommand{\hs}{\hat{\sigma}^2_{(s)}}
\newcommand{\hu}{\mathfrak{s}^2_{(s)}}
\newcommand{\ns}{n^{(s)}}
\newcommand{\Ns}{N^{(s)}}
\newcommand{\mx}[1]{\bar{X}^{(s)}_{#1\cdot}}
\newcommand{\cd}{\ \stackrel{{\bf d}}{\longrightarrow\ }}
\newcommand{\cp}{\ \stackrel{{\bf P}}{\longrightarrow\ }}
\newcommand{\cas}{\ \stackrel{{\bf a.s.}}{\longrightarrow\ }}
\newtheorem{theorem}{Theorem}
\newtheorem{lemma}{Lemma}
\newtheorem{proposition}{Proposition}

\newcommand{\bmu}{\mathbf{\mu}}
\newcommand{\bP}{{\mathbf P}}
\newcommand{\Rbold}{\mathbb{R}}
\newcommand{\bone}{\mathbf{1}}
\newcommand{\noopsort}[1]{}

\begin{document}

%\begin{frontmatter}
\title{Variance Estimation for Tree Order Restricted Models}

%\titlerunning{Variance estimation under tree order restriction}
%\begin{aug}
%\author{\fnms{Antar} \snm{Bandyopadhyay}\thanksref{a}\ead[label=e1]{antar@isid.ac.in}}
%\and
%\author{\fnms{Sanjay} \snm{Chaudhuri}\thanksref{b}\ead[label=e2]{sanjay@stat.nus.edu.sg}\thanks{This work was partially supported by the grant number $R155000111112$ from National University of Singapore.}}

%\address[a]{Indian Statistical Institute, Delhi and Kolkata. \printead{e1}}
%\address[b]{Department of Statistics and Applied Probability, National University of Singapore, Singapore. \printead{e2}}
%\affiliation{Indian Statistical Institute, Delhi and Department of Statistics and Applied Probability, National University of Singapore, Singapore}

%\end{aug}

%\author{Antar Bandyopadhyay $^{\rm a}$ and Sanjay Chaudhuri $^{\rm b}$$^{\ast}$.\thanks{$^{\ast}$ This work was partially supported by the grant number $R155000111112$ from National University of Singapore.}\\\vspace{6pt} $^{\rm a}${\em{Indian Statistical Institute,Delhi and Kolkata.}} Email{antar@isid.ac.in}\\\vspace{6pt}$^{\rm b}${\em{Department of Statistics and Applied Probability, National University of Singapore, Singapore.} Email{sanjay@stat.nus.edu.sg}}}

\author{Antar Bandyopadhyay}
\address{Indian Statistical Institute,Delhi and Kolkata.}  
\email{antar@isid.ac.in} 
\author{Sanjay Chaudhuri}
\address{Department of Statistics and Applied Probability,\\National University of Singapore, Singapore.}
\email{sanjay@stat.nus.edu.sg}
\thanks{This work was partially supported by the grant number $R155000111112$ from National University of Singapore.}

%\and Department of Statistics and Applied Probability, National University of Singapore, Singapore \email{sanjay@stat.nus.edu.sg}} 
%\author{Sanjay Chaudhuri}
%\institute{\at Department of Statistics and Applied Probability,\\National University of Singapore, Singapore.\\sanjay@stat.nus.edu.sg}

%\date{April 27, 2012}
%\begin{document}
\maketitle
\begin{abstract}
In this article we discuss estimation of the common variance of several normal populations with tree order restricted means. 
We discuss the asymptotic properties of the maximum likelihood estimator of the variance as the number of populations tends to infinity.
We consider several cases of various orders of the sample sizes and show that the maximum likelihood estimator of the variance may or may not be consistent or be asymptotically normal.
\end{abstract}

%\end{abstract}
%\begin{keyword}
%\kwd{Order restricted inference}
%\kwd{tree order restriction}
%\kwd{maximum likelihood estimation}
%\kwd{high dimensional data}
%\kwd{variance estimation}
%\end{keyword}
%\end{frontmatter}

%\maketitle

\section{Introduction}
\label{sec:Intro}

\subsection{Background and Motivation}
\label{subsec:motiv}
Tree order restrictions arise naturally in many important applications. One classical scenario is the comparison of several, say $s$, treatments with a known control or a placebo treatment.
It is then natural to model the effect of the $i^{\mbox{th}}$ treatment, say $\mu_i$ to be at least as large as the effect of the control treatment denoted by $\mu_0$, that is, $\mu_0\le\mu_i$, for all $i=1,2,\ldots,s$.
 
Under the tree order restriction the parameter space for 
$\bmu := \left(\mu_0; \mu_1,\right.$ $\left.\mu_2, \ldots, \mu_s\right)$ 
forms a symmetric polyhedral cone in $\Rbold^{s+1}$
with its spine along the line $\mu_0=\mu_1=\cdots=\mu_s$. It can be shown that under the normality assumption the constrained \emph{maximum likelihood estimator} (MLE) of the mean vector 
$\bmu$ is biased in many situations.  
In fact in \citet{leecic1} it was shown that if $\mu_i$'s and the sample sizes from each population remain bounded then 
the bias for $\mu_0$ diverges to $-\infty$ as $s\rightarrow \infty$. 
Because of this phenomenon the constrained MLEs have been criticised severely in literature. \citet{hwang_peddada_1994} wrote that the MLE ``fails disastrously'' and  \citet{cohen_sackrowitz_2002} remarked that the MLE is ``undesirable''. 
%Some estimators for order restricted means have been studied by \citet{vaneeden_zidek_2002}.
  
Under certain conditions however, $\mu_0$ is unbiased. It was shown by 
\citet{chaudhuri_perlman_spl_2005,chaudhuri_perlman_jspi_2007}
that if either $\mu_{(1)} :=\min_{i\ge 1}~\mu_i$ or the sample size from the population corresponding to $\mu_0$ grow sufficiently fast, then the
MLE of $\mu_0$ can be bounded from below in probability and it may even be consistent. 

In all earlier works
it is generally assumed that the treatment groups are homoskedastic but none considered estimation of the variance $\sigma^2$.  
In this article we discuss maximum likelihood estimation of $\sigma^2$ under normality assumption with the tree order restriction on the population means.   
We consider the asymptotic properties of this estimator as $s\rightarrow \infty$,  that is, in the limit the dimension of $\bmu$ becomes large.  The sample size drawn from each population (denoted $\ns_i$ for the $ith$ population) is assumed to grow with $s$ at various rates. 

The main findings of this article can be summarised as follows. First of all, we show a curious phenomenon that depending on the growth of $\ns_i$ with $s$ the MLE $\hs$ is consistent under mild conditions, even though
in some of these cases $\hm$ is not consistent. 
Under stricter assumptions we also prove asymptotic normality for the estimator $\hs$. We show that the so called Neyman-Scott phenomenon extends to our setup.  We perform simulation studies for the cases left unresolved by our assumptions.
Finally, the MLE is compared with the unbiased estimator of variance which ignores the tree order restriction on the mean vector. The simulation study shows that the MLE dominates the latter in terms of the variance. 
However, because of its bias for finite values of $s$, the MLE has a higher mean squared error (MSE) than this unbiased estimator.    

\subsection{Constrained Maximum Likelihood Estimator of $\mu$ and $\sigma^2$}
\label{subsec:MLEs}
Suppose there are $s+1$ independent Normal populations indexed by $0$, $1$, $2$, $\ldots$, $s$ with unknown means 
$(\mu_i)_{i\ge 0}$ and unknown common variance $\sigma^2$.   
Let $X_{i1}$, $X_{i2}$, $\ldots$, $X_{i\ns_i}$, be an i.i.d. sample of size $\ns_i$ from the $i^{\mbox{th}}$ population.  Further suppose that $\Ns=\sum^s_{i=0}\ns_i$ be the total sample size and we write $\mx{i}:=\sum^{\ns_i}_{j=1}X_{ij}/\ns_i$ for the sample mean of the $i^{\mbox{th}}$ population where
$ 1 \leq i \leq s$.

Under the tree order restriction and the assumptions made above, the maximum likelihood estimator of $\mu$ is given by 
(see \citep{robertson_wright_dykstra_book, silvapulle_sen_book})
\begin{align}
\hm_0&=\min_{I\subseteq\{1,2,\ldots,s\}}\frac{\ns_0\mx{0}+ \sum_{i\in I} \ns_i \mx{i}}{\ns_0+\sum_{i \in I} \ns_i} \,\,\, \mbox{and} \label{eq:mu0}\\
\hm_i&=\max\left(\hm_0,\mx{i}\right) \,.\label{eq:mui}
\end{align}
Note that, if we do not assume Normality even then the
equations \eqref{eq:mu0} and \eqref{eq:mui} give the \emph{constrained least squared estimators} of the
mean vector.

It is then immediate that
under the tree order restriction and Normality assumption the constrained maximum likelihood estimator of $\sigma^2$ is given by
\begin{equation}\label{eq:lglkhd}
\hs=\frac{1}{\Ns}\left\{\sum^{\ns_0}_{j=1}\left(X_{0j}-\hm_0\right)^2+\sum^s_{i=1}\sum^{\ns_i}_{j=1}\left(X_{ij}-\hm_i\right)^2\right\}.
\end{equation}
Once again even without the assumption of Normality this also gives the constrained least squared estimator of $\sigma^2$. 
%It is obvious that $\hs$ depends on the constrained MLE of $\bmu$.
 
\subsection{Main Results}
\label{subsec:results}
We make the following assumptions throughout this article:
\begin{itemize}
\item[{\bf A1:}] The mean vector $\bmu$ is tree order restriction, that is, $\mu_0\le\mu_i$ for all $i\geq 1$.
\item[{\bf A2:}] There exists $B > 0$ but \emph{unknown} such that $\mu_i \leq B$ for all $i \geq 0$. 
\item[{\bf A3:}] The populations are Normal. 
\end{itemize}
For simplicity we further assume that
\begin{itemize}
\item[{\bf A4:}] $\ns_1=\ns_2=\cdots=\ns_s=\ns$ and both $\ns_0$ and $\ns$ are non-decreasing in $s$.
\end{itemize}
We would like to note here that our assumption {\bf A2} does not put any further constraint on the parameter space for computation of the MLEs
since we assume $B$ is unknown. This is quite different from what was assumed by \citet{barlow_bartholomew_1972} 
for bounded isotonic regression model. They assumed the bound to be known. 

Our main interest is to study the asymptotic properties of the maximum likelihood estimator of $\sigma^2$ namely $\hs$ as the number of populations becomes 
large.

%Note that, $\mu_i\le B$, for all $i\ge 0$ is a part of the condition given in \citep{leecic1} for the divergence of the MLE of $\mu_0$.  
%Under this condition the MLE of $\mu_0$ diverges, 
%unless $\ns_0$ grows sufficiently fast with $s$ \citep{chaudhuri_perlman_jspi_2007}.  

We first consider an example with two populations, with tree-ordered means. The size of the sample drawn from the population with larger mean increases linearly with $s$, while size of the placebo sample remains constant. 
  
\begin{theorem}
\label{lem:fnts}
Consider two populations with a common variance $\sigma^2$ and means $\mu_0\le\mu_1$. 
Let $\ns_0=m$ and $\ns_1=m^{\prime}s$, where $m$, $m^{\prime} \geq 1$ are two fixed integers. Then as $s\rightarrow\infty$
\begin{equation}
\hs \longrightarrow \sigma^2 \,\,\, \mbox{a.s.}
\label{Equ:Consistency}
\end{equation}
Moreover,
\begin{equation}
\sqrt{\Ns} \left(\hs - \sigma^2 \right) \cd \mbox{N}\left(0, 2\sigma^4\right) \,. 
\label{Equ:CLT}
\end{equation}
\end{theorem}
Note that as pointed out in Section \ref{sec:proofs}
it is not difficult to show that $\hm_0$ is biased in this case, yet $\hs$ is consistent and a \emph{central limit theorem} (CLT) holds. 

The assumptions of Theorem \ref{lem:fnts} can be interpreted in the following alternative way.  Suppose we consider $s+1$ populations with an unknown common variance $\sigma^2$ and $\mu_1=\mu_2=\cdots=\mu_s$ with $\mu_0\le\mu_1$. Both $\mu_0$ and $\mu_1$ are unknown.  
Let $\ns_0=m$ and $\ns=m^{\prime}$ be the sample sizes from these distributions.  Clearly the MLEs of $\mu_0$, $\mu_1$ and $\sigma^2$ are exactly same as in Theorem \ref{lem:fnts}.  
So it follows that in the limit as $s\rightarrow \infty$, the MLE $\hs$ is consistent and admits a CLT. It is worth mentioning here that 
the assumption that $\mu_1=\mu_2=\cdots=\mu_s$ is very crucial for the consistency and also for the CLT. This is because as stated in Theorem \ref{thm:asf} below the consistency may fail and the simulations
presented in Section \ref{sec:simul} shows that CLT may not hold either. 

Our next theorem deals with the case when
the total sample size $\Ns$ grows at a faster rate than $s$.
\begin{theorem} 
\label{thm:general}
Suppose $\Ns \rightarrow \infty$ then under the assumptions {\bf A1 -- A4}
\begin{equation}
\bP\left(0 \leq \limsup_{s \rightarrow \infty} \hs \leq \sigma^2\right) = 1 \,.
\label{Equ:General-LimSup-Bound}
\end{equation}
Further if we assume that
$s/\Ns \longrightarrow 0$ as $s \rightarrow \infty$ then 
\[
\hs \rightarrow \sigma^2 \,\,\, \mbox{a.s.} \,,
\]
while if $s/\sqrt{\Ns} \longrightarrow 0$ then
\[
\sqrt{\Ns}\left(\hs - \sigma^2\right) \cd \mbox{N}\left(0, 2\sigma^4\right) \,.
\]
\end{theorem}
From the proof of Theorem \ref{thm:general} in Section \ref{sec:proofs} we see that the result
holds for any constraint on the mean vector $\bmu$ which need not be just the tree order restriction.
So this result may be used in other constrained problems, for example, in the study of the
\emph{isotonic regression model}, where one assumes that
$\mu_0 \leq \mu_1 \leq \mu_2 \leq \cdots \leq \mu_s$.

Following result is immediate from Theorem \ref{thm:general} which
covers the case when $\ns\equiv m$ fixed but $\ns_0$ is increasing at an appropriate rate. 
\begin{theorem}
\label{thm:hmcons}
Let $\ns_0$ be such that $s/\ns_0 \longrightarrow 0$ as $s \rightarrow \infty$. Then under the assumptions {\bf A1 -- A4} $\hs$ is strongly consistent. 
Moreover, if $s/\sqrt{\ns_0} \longrightarrow 0$ as $s\rightarrow\infty$ then,
\[
\sqrt{\Ns}\left(\hs-\sigma^2\right)\cd \mbox{N}\left(0, 2\sigma^4\right).
\]
\end{theorem}
It is worth noting that under the conditions in Theorem \ref{thm:hmcons} the MLE $\hm_0$ is consistent \citep{chaudhuri_perlman_jspi_2007}.

From Theorem \ref{thm:general} the strong consistency holds if either $\ns_0/s \longrightarrow \infty$ or $\ns \longrightarrow \infty$ as $s \rightarrow \infty$. For the CLT to hold we need stronger condition, namely, 
$\ns_0/s^2 \longrightarrow \infty$ or
$\ns/s \longrightarrow \infty$ as $s \rightarrow \infty$. 
In particular it covers the case when $\ns/\log s \longrightarrow \infty$ for which 
$\hm_0$ is consistent if and only if $\ns_0 \rightarrow \infty$ 
(see Proposition \ref{prop:iff} in Section \ref{sec:technical}).  In this case we have not been able to proof the
CLT, which may hold (see Section \ref{sec:simul} for simulated results). 

Following theorem deals with the case when 
both $\ns_0$ and $\ns$ remain bounded. 
\begin{theorem}\label{thm:asf}
Suppose $\ns_0=\ns=m$ for some fixed $m\ge 2$. Then under the assumptions {\bf A1 -- A4}
\begin{equation}
\hs \longrightarrow \frac{m-1}{m}\sigma^2 \,\,\, \mbox{a.s.}
\label{Equ:Inconsistency}
\end{equation}
\end{theorem}
Like in Theorem \ref{lem:fnts} in this case also $\Ns \sim m s$  but here $\hs$ is not consistent.   
The difference is in the dimension of the mean parameter being estimated. 
In Theorem \ref{lem:fnts} it is exactly $2$, however in Theorem \ref{thm:asf} it grows unbounded with $s$. 
It is also worth noting that in this case $\hm$ is infinitely biased, that is, $\hm \longrightarrow \infty$ a.s. \citep{leecic1, chaudhuri_perlman_jspi_2007}
yet the asymptotic bias of $\hs$ is small if $m$ is large.
Once again in this case the CLT may not hold and we present some simulation results in
Section \ref{sec:simul}.

The apparent ambiguity between Theorems \ref{lem:fnts} and \ref{thm:asf} is reminiscent of the so called Neyman-Scott example \citep{neyman_scott_1948}.  
They considered i.i.d. samples of equal finite size from several normal populations with unknown means and common variance.  
It was shown that in the limit if number of populations increases the MLE of the common variance is inconsistent.  Here we observe the same 
phenomenon with tree order restriction on $\mu$.  For estimation of $\sigma^2$, $\mu$ is a nuisance parameter.  
In Theorem \ref{lem:fnts}, $\mu_i=\mu_1$, for all $i\ge 1$.  Thus even in the limit of $s\rightarrow \infty$, the number of nuisance parameters remain 
bounded.  
In contrast, in Theorem \ref{thm:asf} this number increases unbounded. This explains the inconsistency of 
$\hs$ in the latter.  
\citet{chaudhuri_perlman_jspi_2005} argue that in general, if the number of nuisance parameters is allowed to grow, the MLE of the parameter of interest
may not be consistent. It even may not converge to any limit.  

It is also worth noting that our model does not fall into the general class of model discussed in \citep{murphy_vaart_2000}. This is because in our case
the total sample size, namely, $\Ns$ depends on the number of populations $s$ and we consider asymptotics as $s \rightarrow \infty$. 
Thus we observe some non-standard limiting results unlike the case in \citep{murphy_vaart_2000}.

One important case which is not covered by the results described above is when $\ns_0=O(s)$ and $\ns$ remains bounded. 
\citet{chaudhuri_perlman_jspi_2007} show that in this case $\hm_0$ remains bounded from below with high probability, but may not be consistent.  
In Section \ref{sec:simul} we present some simulation results for this case.

Before we end this subsection we would like to point out that from the proof of Theorem \ref{lem:fnts} it is clear that the 
Normality assumption {\bf A3} is not needed for this result. Moreover the general result Theorem \ref{thm:general} 
and hence its corollary Theorem \ref{thm:hmcons}
can also be derived without using assumption {\bf A3}, but just assuming that the data comes for a \emph{location-scale family} with mean
$\mu_i$ for the $i^{\mbox{th}}$ population and common variance $\sigma^2$ and finite forth moment. But this only complicates the proof
and thus we present a proof under assumption {\bf A3}. It is not clear though whether the statement of Theorem \ref{thm:asf} follows
without assumption of the Normality for the populations (assumption {\bf A3}). It is only in this proof we truly need this assumption. 

\subsection{Comparison with an unconstrained estimator of $\sigma^2$}\label{sec:ub}
If the information about the tree order restriction on $\mu$ is ignored, a natural estimator of $\sigma^2$ is given by (see \citet{peddada_et_al_2006}):
\begin{equation}
\hu=\frac{1}{\Ns-(s+1)}\sum^s_{i=0}\sum^{\ns_i}_{j=1}\left(X_{ij}-\mx{i}\right)^2.
\end{equation}  
It is well known that under mild conditions eg. existence of second moments $\hu$ is unbiased and strongly consistent for $\sigma^2$.
%Under our normality assumption it trivially follows that $\{\Ns-(s+1)\}\hu/\sigma^2\sim\chi^2_{(\Ns-s-1)}$ distribution, $E[\hu]=\sigma^2$ and $Var[\hu]=2\sigma^4/(\Ns-s-1)$. 
\begin{theorem}\label{thm:hu} The following results hold:
\begin{enumerate}
\item Under assumptions {\bf A1 -- A4}, $E[\hs]\le E[\hu]$.
\item $\hs\ge (\Ns-s-1)\hu/\Ns$ always holds.
\item Under the assumptions of Theorem \ref{lem:fnts}, Theorem \ref{thm:general}  when $s/\Ns\rightarrow 0$ and \ref{thm:hmcons},
\[\mid\hs-\hu\mid\cas 0.
\]
\end{enumerate} 
\end{theorem}

Theorem \ref{thm:hu} shows that $\hs$ is biased and may be smaller than the $\hu$ for finite $s$, however, when $\Ns$ grows faster than the dimension, $\hu-\hs$ converges in probability to $0$.  

It should be noted that, no analytic expression of bias and variance of $\hs$ is available for finite $s$. Thus the mean squared error (MSE) of $\hs$ and $\hu$ cannot be compared analytically. Furthermore, the case when $\Ns\sim ms$ (Theorem \ref{thm:asf}) remains unresolved in Theorem \ref{thm:hu}.  
These questions are explored in a simulation study presented in the next section.
\subsection{Outline} 
\label{subsec:outline}
The rest of the article is structured as follows. 
The next section gives the detailed simulation results of the cases mentioned above.
In Section \ref{sec:proofs} we present the proofs of the main results.
Section \ref{sec:technical} contains some of the technical results and their proofs which we 
use to prove the main results.

\section{Simulation Studies for Some Unresolved Cases}\label{sec:simul}

\begin{figure}[t]
\begin{center}
\subfigure[\label{fig:1bp}]{
\resizebox{2in}{2in}{\includegraphics{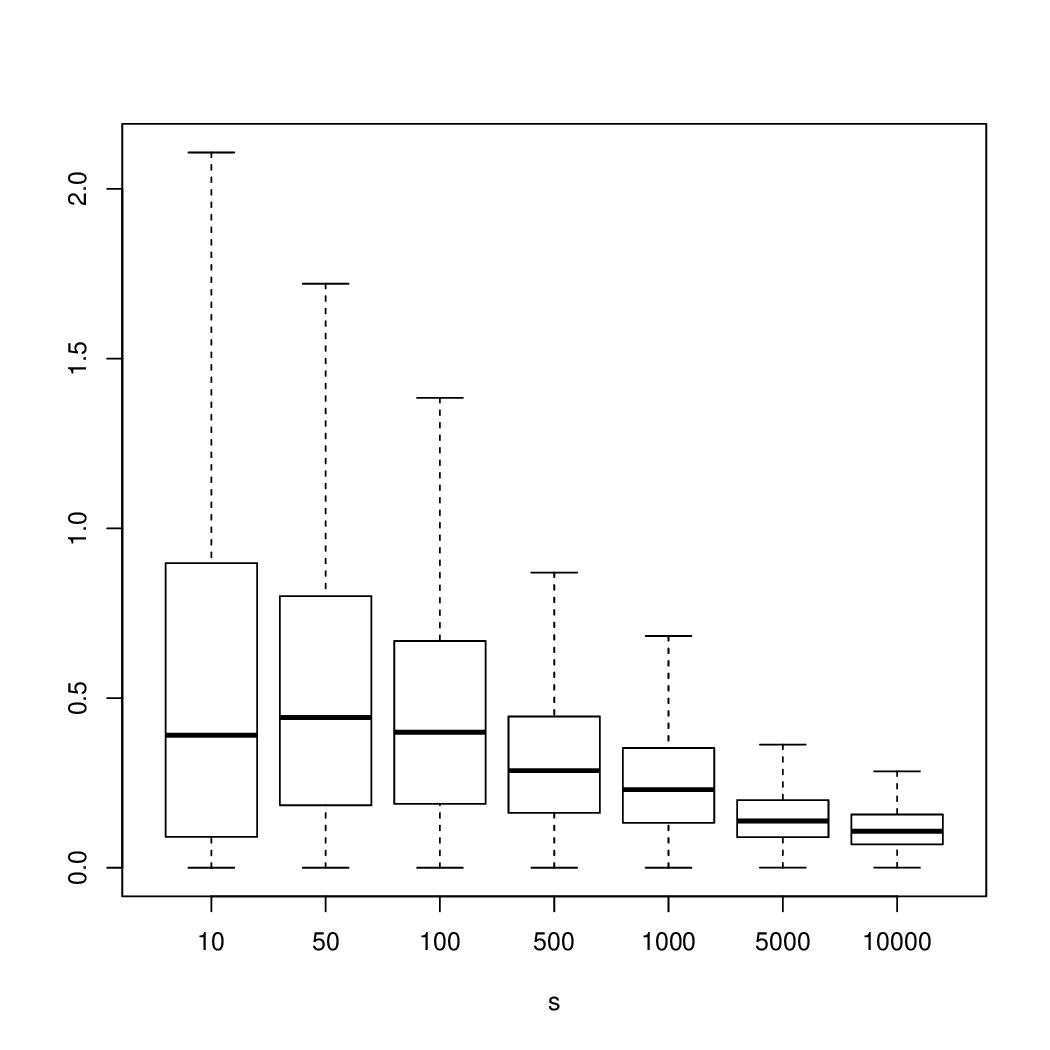}}}\qquad
\subfigure[\label{fig:1hg}]{
\resizebox{2in}{2in}{\includegraphics{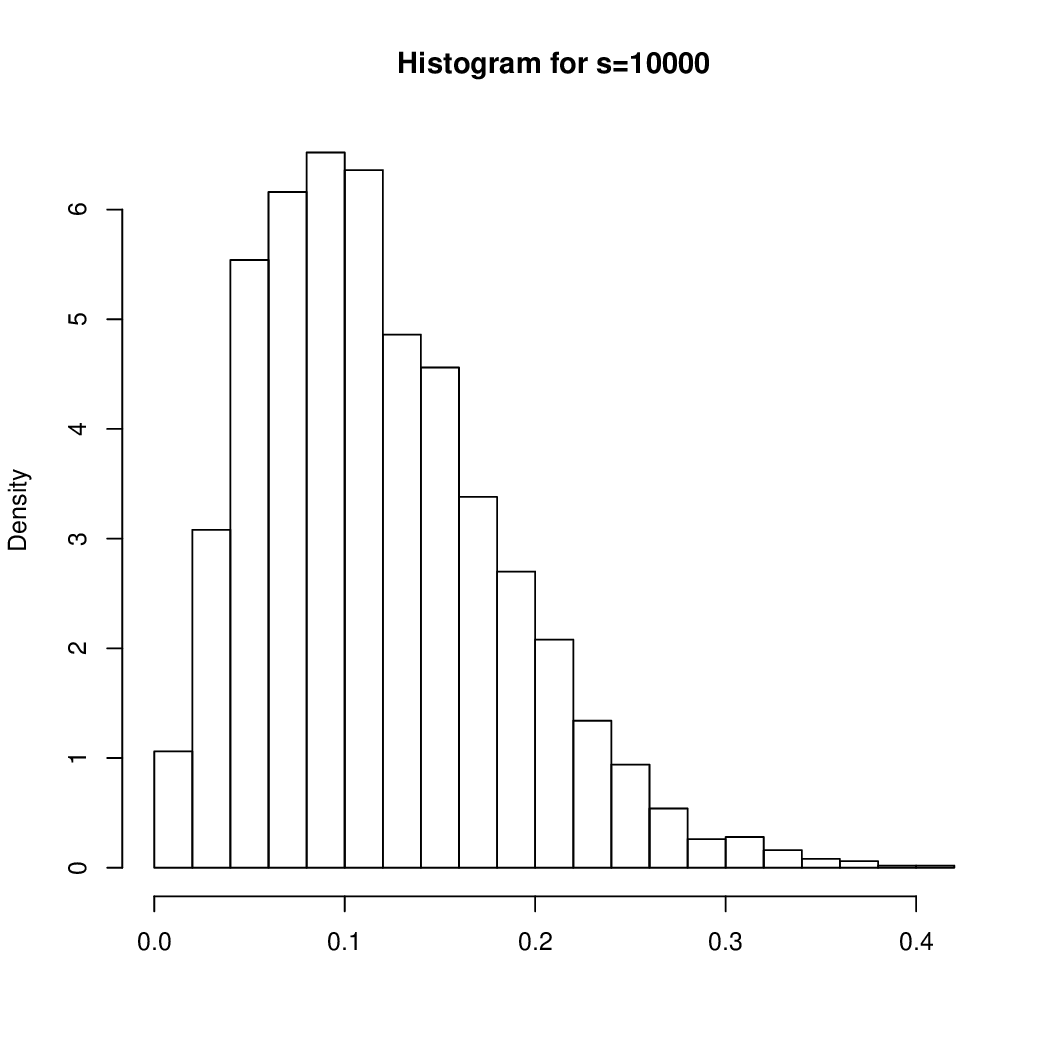}}}
\end{center}
\caption{Box plot (Figure \ref{fig:1bp}) and histogram for $s=10000$ (Figure \ref{fig:1hg}) of $\xi_s=\sum^s_{i=0}\ns_i\left(\mx{i}-\hm_i\right)^2/\sqrt{\Ns}$ when $\ns_0=\ns=m$.}
\label{fig:1}
\end{figure}

In this section we perform a simulation study to explore two sets of properties of $\hs$.  First is its asymptotic behaviour under some conditions which are not covered by the results in Section \ref{subsec:results}.  Second, we consider its bias and variance for finite $s$ and compare them with that of $\hu$ described in Section \ref{sec:ub}.

The following three situations are not covered by the results in Section \ref{subsec:results}.
\begin{enumerate}
\item We check if the CLT holds when $\ns_0=\ns=m$. Theorem \ref{thm:asf} we know that $\hs$ is inconsistent for $\sigma^2$ in this case. 
\item We check if CLT holds when $\ns_0=\ns=(\log s)^2$. Notice that in this case, $\Ns=(s+1)(\log s)^2$ and from Theorem \ref{thm:general} it follows that $\hs$ is consistent for $\sigma^2$ in this case.  However, $s/\sqrt{\Ns}\not\rightarrow 0$, so the CLT cannot be derived from Theorem \ref{thm:general}. 
\item Asymptotic behaviour of $\hs$, when $\ns_0=O(s)$ and $\ns$ remains bounded is not covered by any result considered above.  In this case even the consistency of $\hm_0$ is not known, though it is bounded below with high probability \citep{chaudhuri_perlman_jspi_2007}. 
\end{enumerate}

\begin{figure}[t]
\begin{center}
\subfigure[\label{fig:2bp}]{
\resizebox{2in}{2in}{\includegraphics{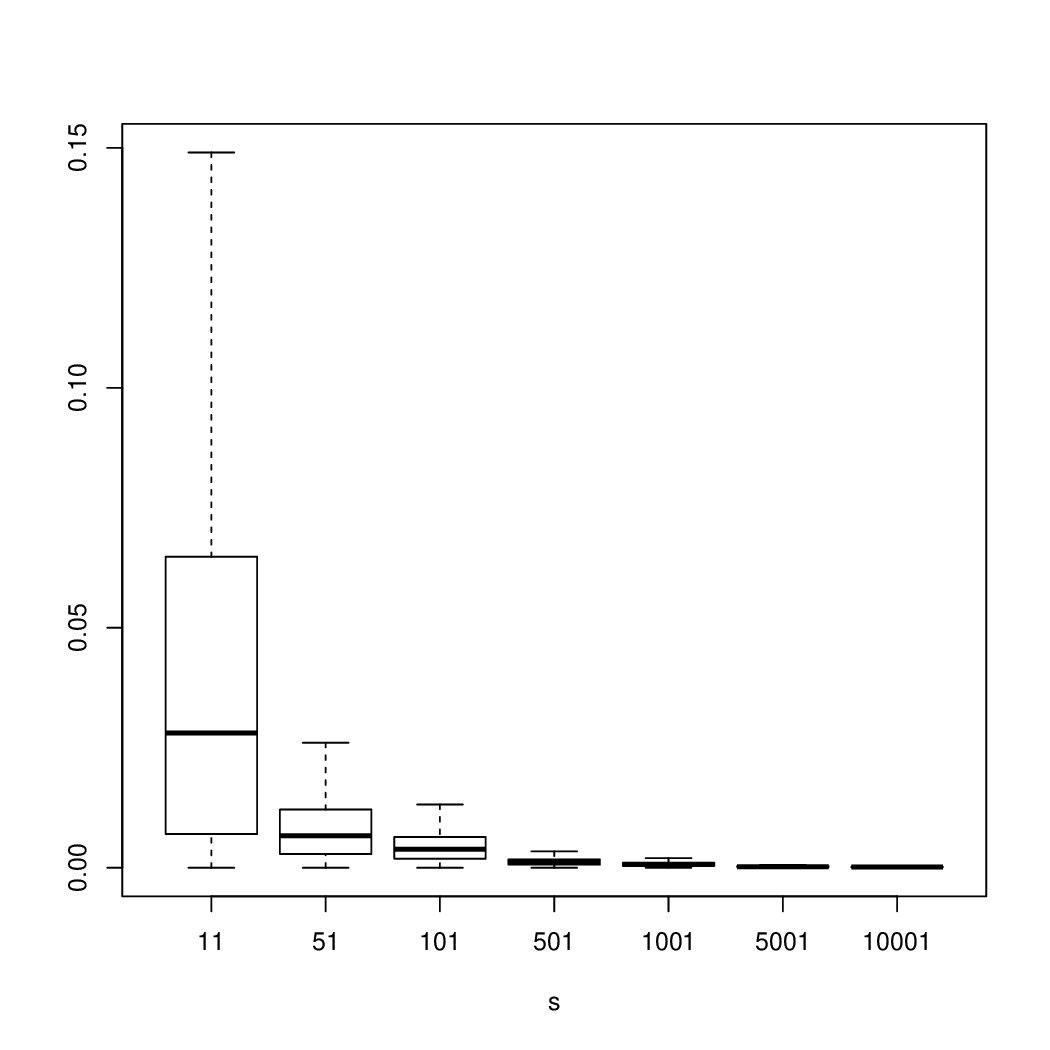}}}\qquad
\subfigure[\label{fig:2hg}]{
\resizebox{2in}{2in}{\includegraphics{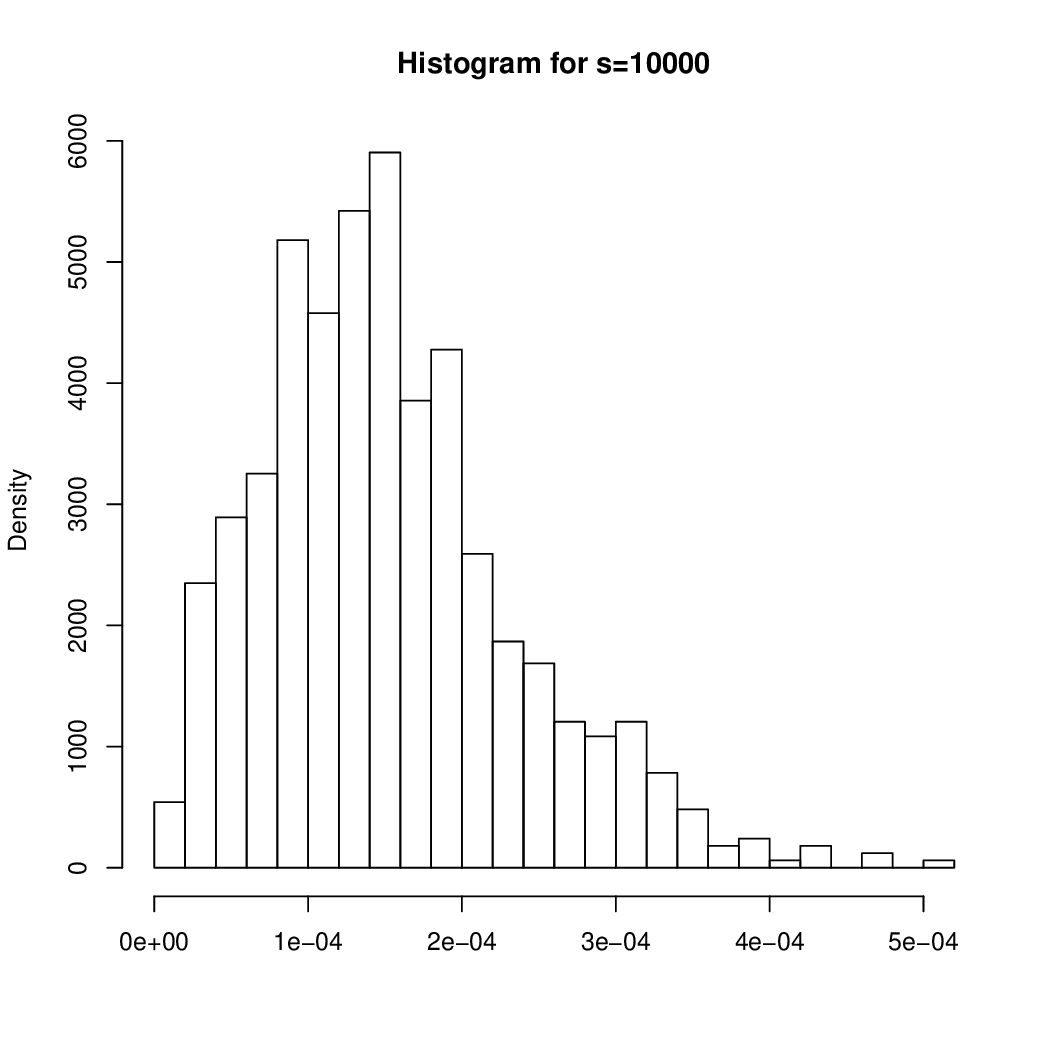}}}
\end{center}
\caption{Box plot (Figure \ref{fig:2bp}) and histogram for $s=10000$ (Figure \ref{fig:2hg}) of $\xi_s$ when $\ns_0=\ns=(\log s)^2$.}
\label{fig:2}
\end{figure}

\begin{figure}[h]
\begin{center}
\resizebox{2in}{2in}{\includegraphics{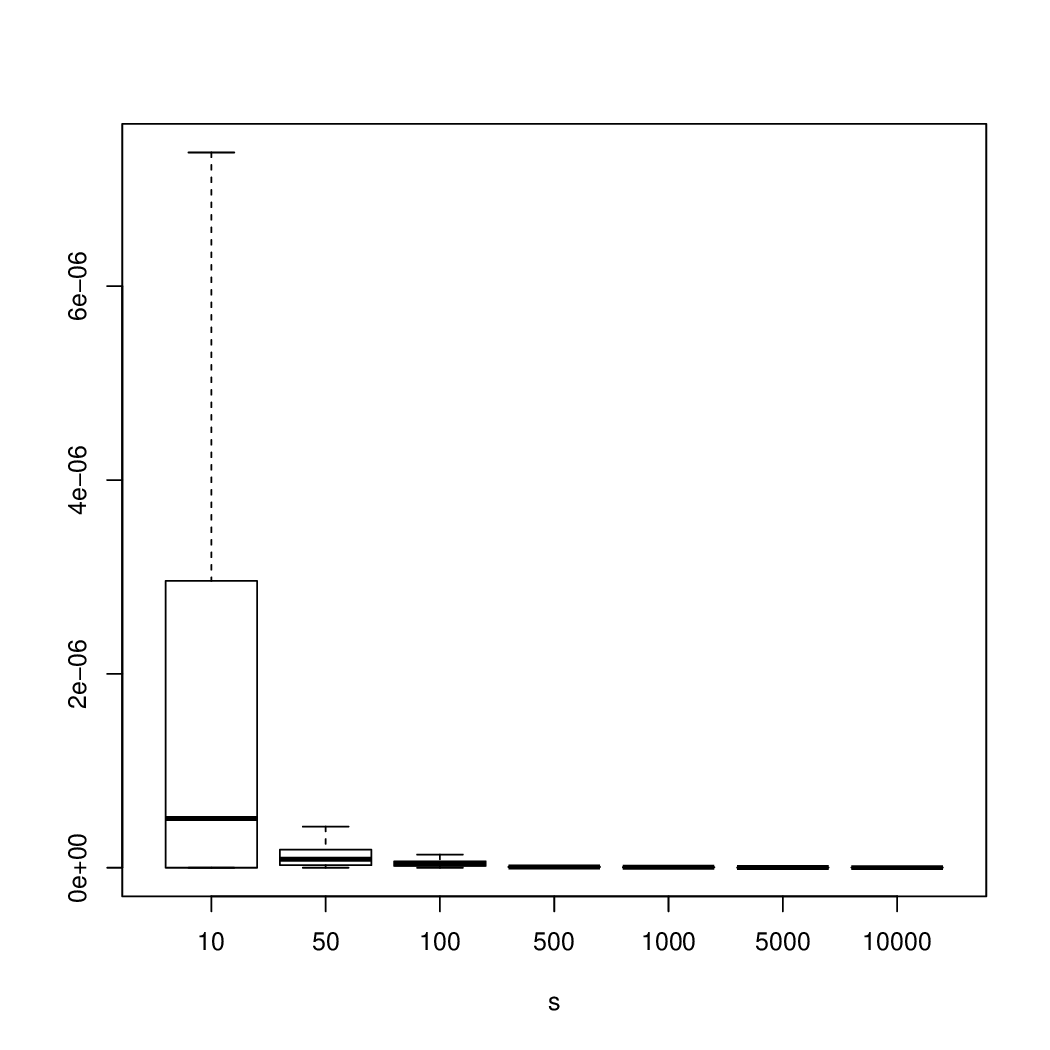}}
\end{center}
\caption{Box plot of $\xi_s/\sqrt{\Ns}=I_2+I_4$ when $\ns_0=O(s)$ and $\ns=m$.}
\label{fig:3u}
\end{figure}

\begin{figure}[h]
\begin{center}
\subfigure[\label{fig:3bp}]{
\resizebox{2in}{2in}{\includegraphics{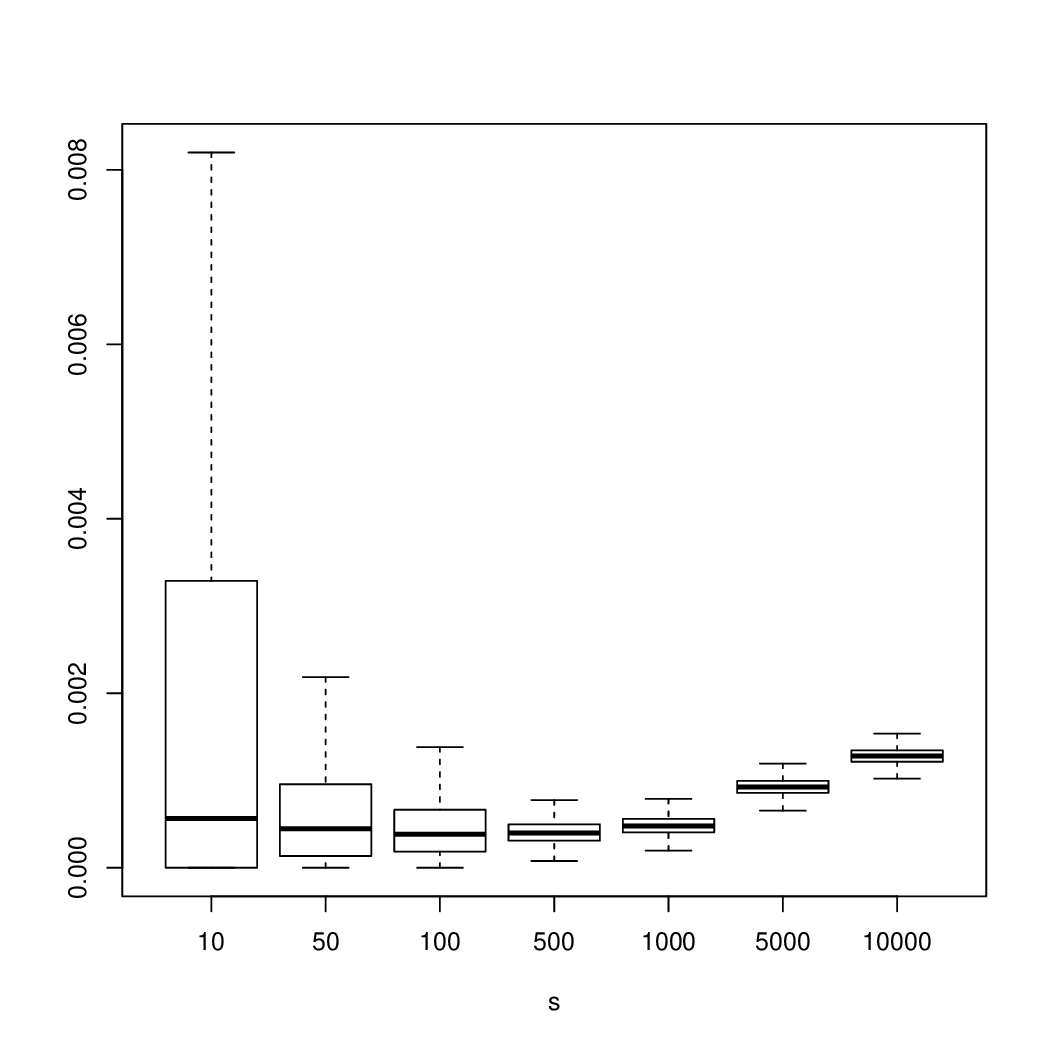}}}\qquad
\subfigure[\label{fig:3hg}]{
\resizebox{2in}{2in}{\includegraphics{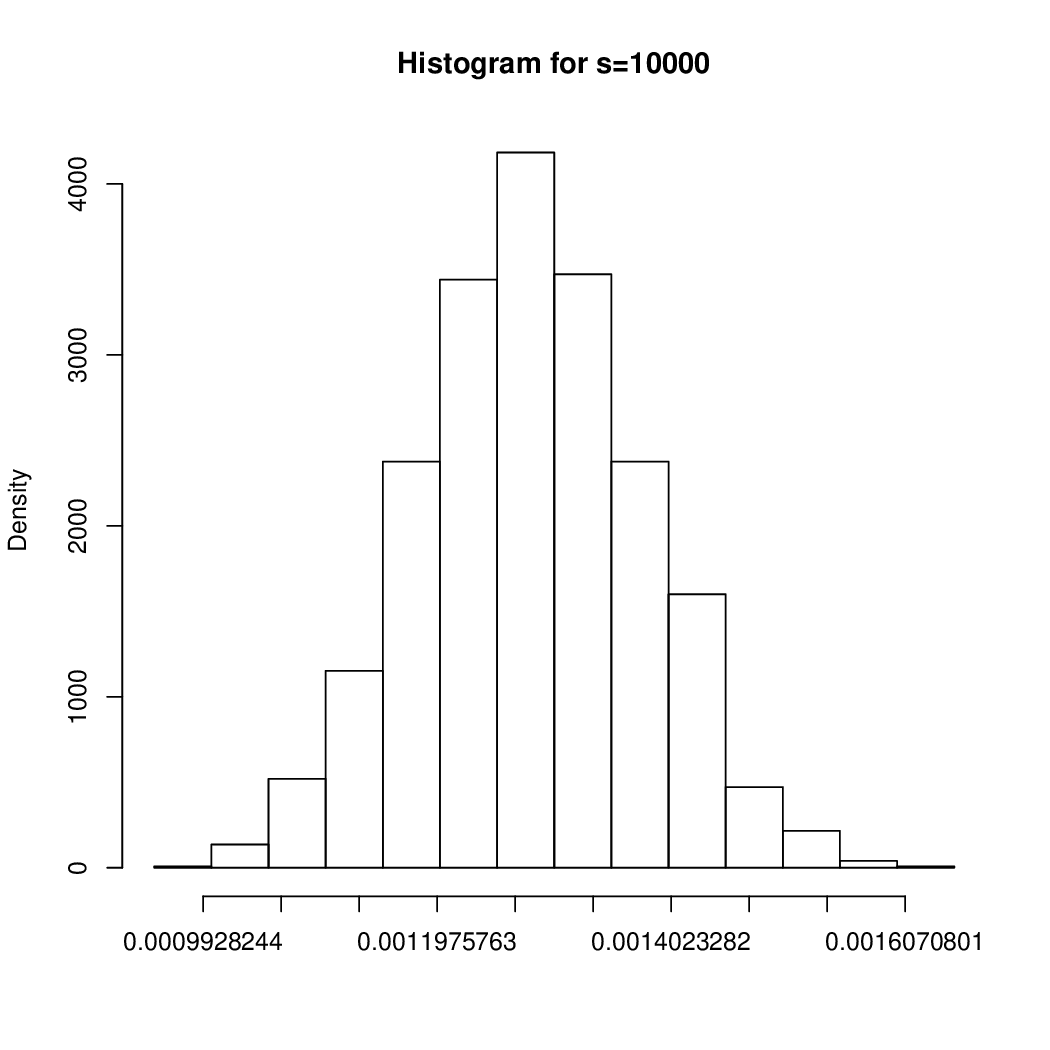}}}
\end{center}
\caption{Box plot (Figure \ref{fig:3bp}) and histogram for $s=10000$ (Figure \ref{fig:3hg}) of $\xi_s$ when $\ns_0=O(s)$ and $\ns=m$}  
\label{fig:4}
\end{figure}

In the simulation study for simplicity we assume all $\mu_i$, $i\ge 0$ to be equal, which is equivalent to assuming $\mu=(0,0,\ldots,0)$. We assume $\sigma^2=1$.  In order to study the asymptotic behaviour for large values of $s$, we consider $s=10$, $50$, $100$, $500$, $1000$, $5000$, $10000$.  The presented results are based on $2500$ repetitions for each population size $s$.

Since, under all the above conditions the term $I_1+I_3$ in \eqref{eq:main2} has a finite limit (in probability) and $\sqrt{\Ns}(I_1+I_3)$ converges in distribution, we concentrate on the random variable $\xi_s=\sqrt{\Ns}(I_2+I_4)=\sum^s_{i=0}\ns_i\left(\mx{i}-\hm_i\right)^2/\sqrt{\Ns}$. 

In the figures we present box and whisker plots of $\xi_s$ for each case. We also present the  histogram of $\xi_s$ for $s=10000$ in all cases.  The top and the bottom of the boxes in the box whisker plots represent the first and third quartiles respectively.  
The line inside the box is the median. The whiskers were drawn to represent the most extreme data point still within the $1.5 $ times the interquartile range of the first and the third quartiles. For better representation we have omitted more extreme points.   

For case (1) ie. when $\ns_0=\ns=m$ from Figure \ref{fig:1bp} it is clear that the spread of the distribution reduces with $s$. Further, it is seen that the medians of the distributions have a decreasing trend with $s$. 
However, we cannot conclude that $\xi_s$ converges to $0$ in probability because the histogram in Figure \ref{fig:1hg} does not show any concentration near $0$.  

The case when $\ns_0=\ns=(\log s)^2$ (ie. case (2)) is presented in Figures \ref{fig:2bp} and \ref{fig:2hg}.  From Figure \ref{fig:2bp} it seems that $\xi_s$ converges to $0$ in probability.  The histogram (Figure \ref{fig:2hg}) also seems to be quite concentrated near $0$. So a CLT may hold in this case.

%Case (3) with $\ns_0=O(s)$ and $\ns=100$ shows a different picture. As seen from the box plot in Figure \ref{fig:3u}, $\xi_s/\sqrt{\Ns}=I_2+I_4$ seems to be rapidly converging to $0$ in probability, which indicates $\hs$ may be consistent. 
%However, the box plot in Figure \ref{fig:3bp} and the histogram Figure \ref{fig:3hg} indicates that asymptotically $\xi_s$ may not be converging to $0$, thus the CLT may not hold.  
%It is interesting to note that the histogram of $\xi_s$ in Figure \ref{fig:3hg} is almost symmetric around its mean, for which we do not have an intuitive explanation. 
 
\begin{figure}[t]     
\parbox[t]{1.7in}{
\resizebox{1.695in}{1.695in}{\includegraphics{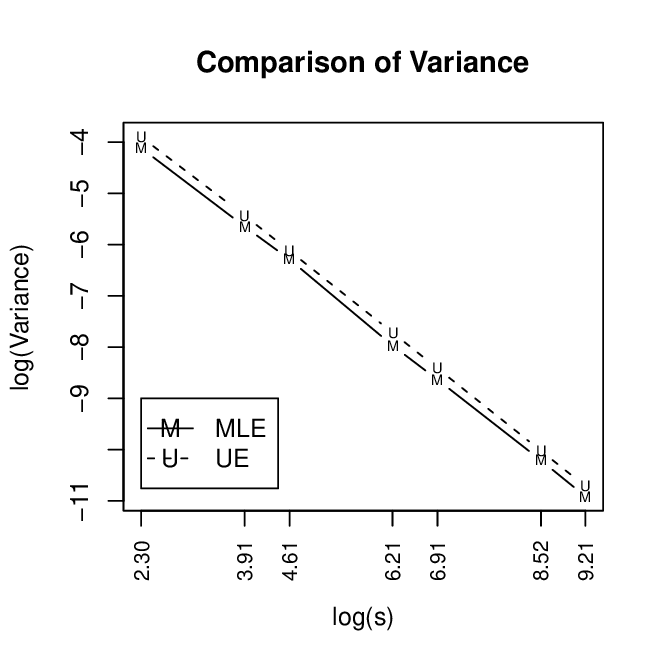}}\\
             \resizebox{1.695in}{1.695in}{\includegraphics{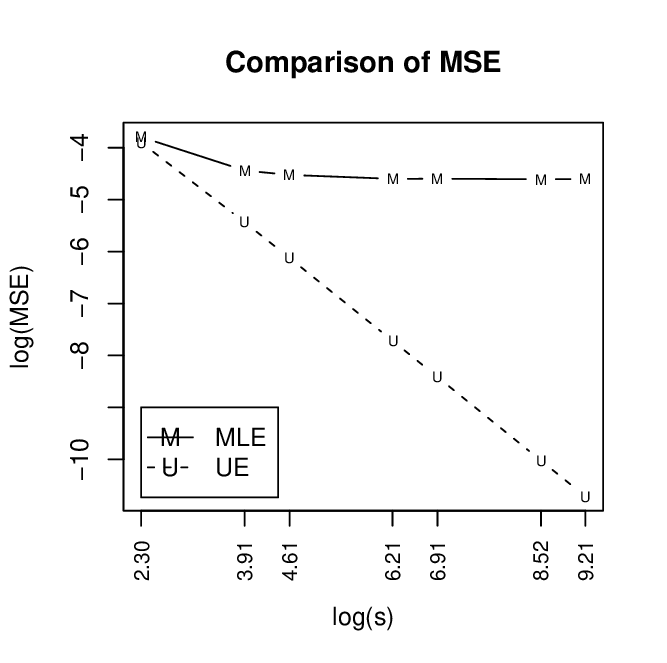}}\\
\centerline{Case $1$}\\
\resizebox{1.695in}{1.695in}{\includegraphics{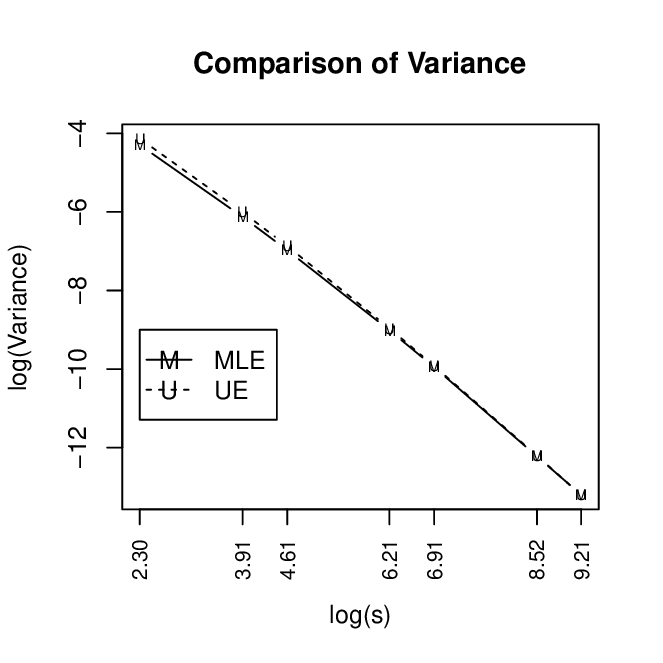}}\\
             \resizebox{1.695in}{1.695in}{\includegraphics{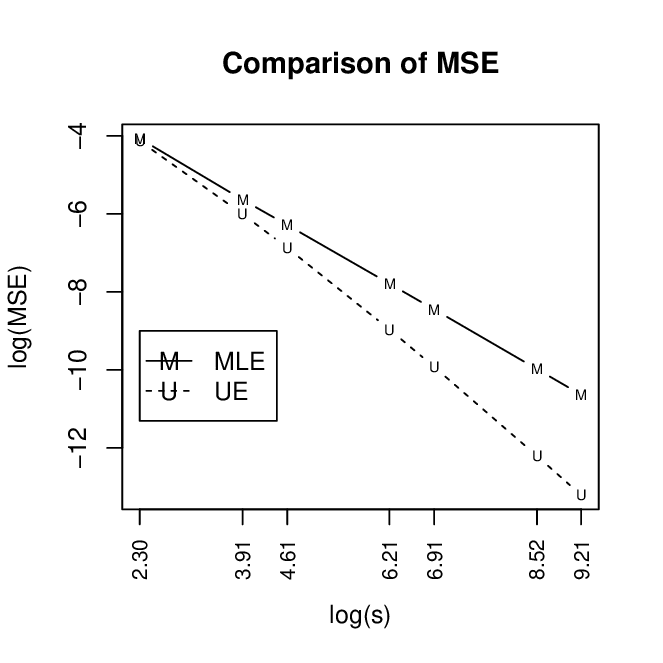}}\\
\centerline{Case $4$}
}\parbox[t]{1.7in}{
\resizebox{1.695in}{1.695in}{\includegraphics{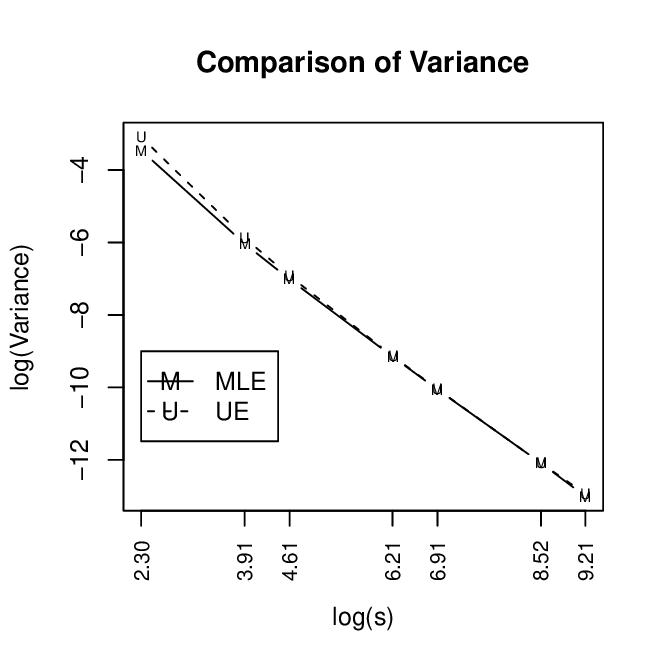}}\\
             \resizebox{1.695in}{1.695in}{\includegraphics{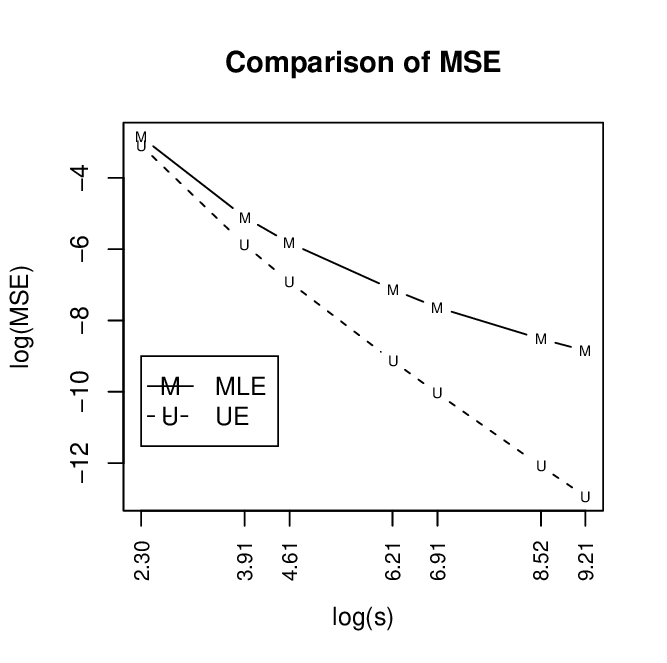}}\\
\centerline{Case $2$}\\
\resizebox{1.695in}{1.695in}{\includegraphics{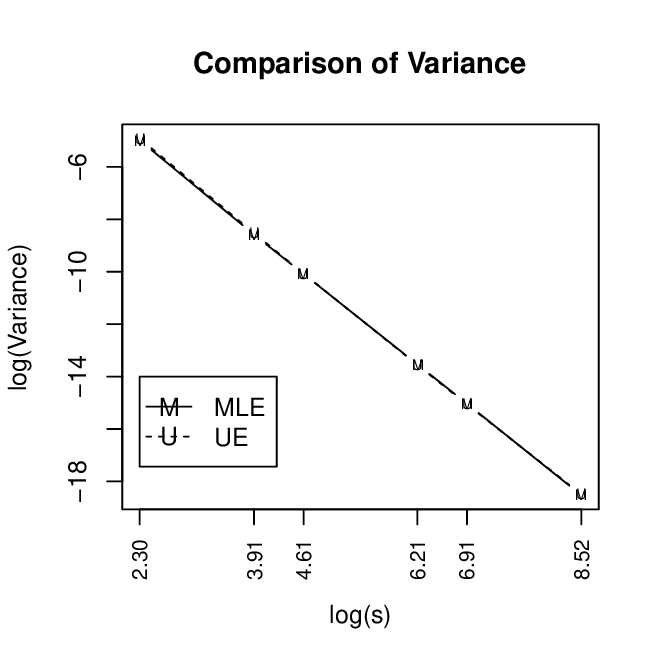}}\\
             \resizebox{1.695in}{1.695in}{\includegraphics{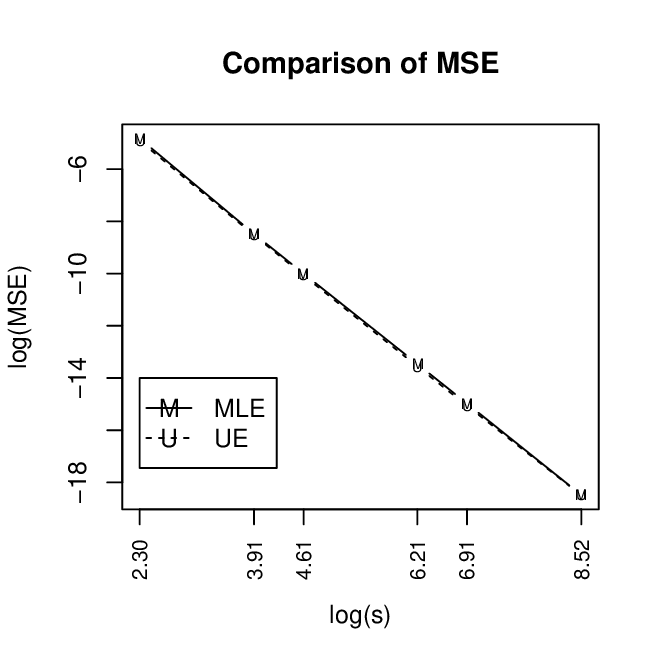}}\\
\centerline{Case $5$}
}\parbox[t]{1.7in}{
\resizebox{1.695in}{1.695in}{\includegraphics{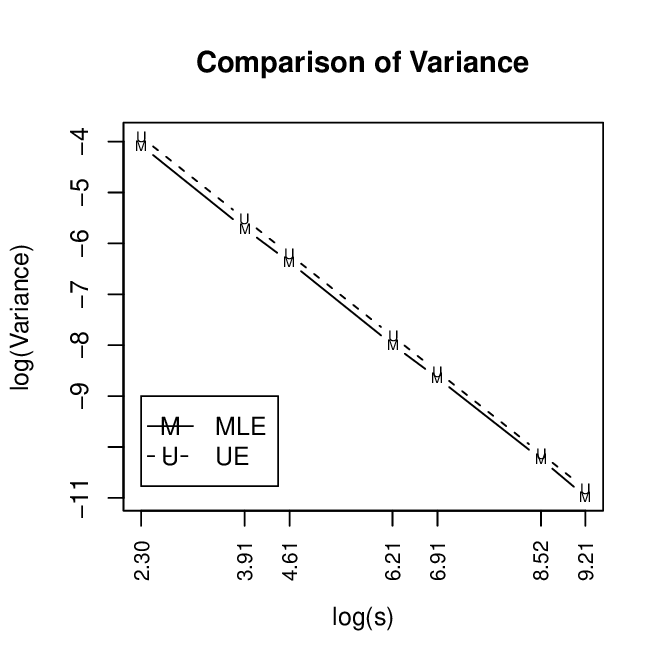}}\\
             \resizebox{1.695in}{1.695in}{\includegraphics{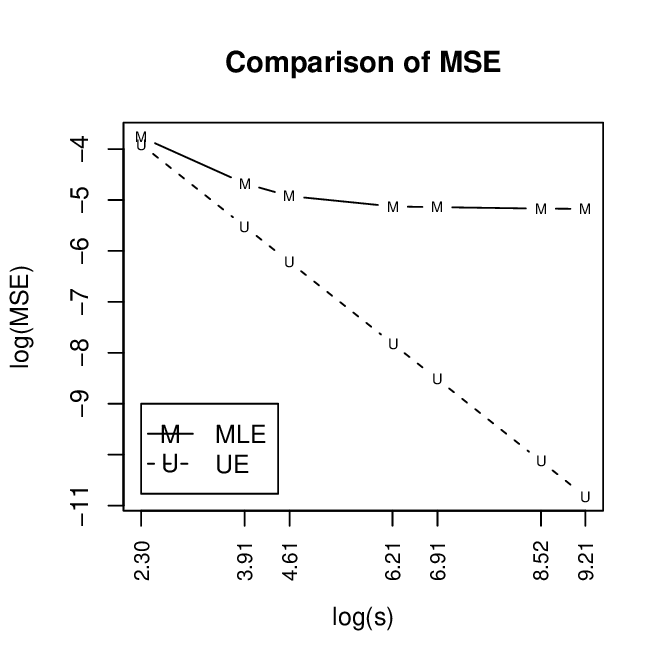}}\\
\centerline{Case $3$}\\
\resizebox{1.695in}{1.695in}{\includegraphics{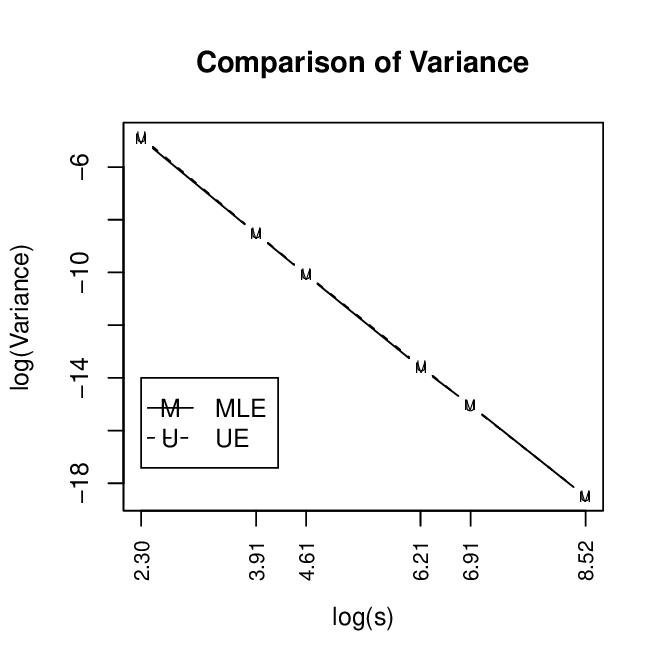}}\\
             \resizebox{1.695in}{1.695in}{\includegraphics{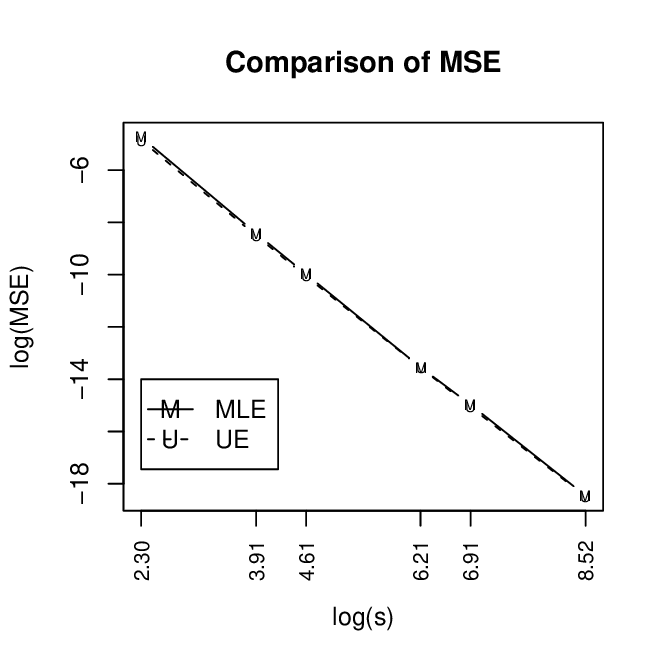}}\\
\centerline{Case $6$}
}
\caption{Comparison of variance and MSE of $\hs$ with $\hu$ in various situations.}
\label{fig:comp}
\end{figure}

Case (3) with $\ns_0=O(s)$ and $\ns=100$ shows a different picture. As seen from the box plot in Figure \ref{fig:3u}, $\xi_s/\sqrt{\Ns}=I_2+I_4$ seems to be rapidly converging to $0$ in probability, which indicates $\hs$ may be consistent. 
However, the box plot in Figure \ref{fig:3bp} and the histogram Figure \ref{fig:3hg} indicates that asymptotically $\xi_s$ may not be converging to $0$, thus the CLT may not hold.  
It is interesting to note that the histogram of $\xi_s$ in Figure \ref{fig:3hg} is almost symmetric around its mean, for which we do not have an intuitive explanation. 
 
In order to compare the variance and the MSE of $\hs$ against $\hu$, in addition to the three cases mentioned above, we considered the following three cases as well:
\begin{enumerate}
\item[$(4)$] $\ns_0=s^{3/2}$, $\ns=10$.  From Theorem \ref{thm:hmcons} it follows that, in this case, $\hs$ is consistent but the CLT is not guaranteed to hold. 
\item[$(5)$] $\ns_0=\ns=s\log(s)$.  In this case, $s/\sqrt{\Ns}\rightarrow 0$.  So by Theorem \ref{thm:general}, $\hs$ is both consistent and a CLT holds. 
\item[$(6)$] $\ns_0=10$, $\ns=s\log(s)$.  Even though $\ns_0$ remains fixed and $\hm_0$ is biased in this case (see Proposition \ref{prop:iff}), $s/\sqrt{\Ns}\rightarrow 0$ and $\hs$ is both consistent and admits a CLT.  
\end{enumerate}

The results are presented in Figure \ref{fig:comp}, where we plot the logarithm of the variance or the MSE against $\log(s)$.  For Cases $5$ and $6$, due to heavy computation, we restrict $s$ up to $5000$.
  In all cases $\hs$ has lower variance than $\hu$. However, due to nonzero bias of $\hs$ for finite $s$, it is dominated by $\hu$ in terms of MSE.  The differences in their MSEs are negligible when CLT holds in Theorem \ref{thm:general} (ie. Case $5$ and Case $6$).  When $\hs$ is biased its MSE naturally tends to the square of the bias (eg. Case $1$).  
Case $3$ shows a very slow downward trend in MSE, which may indicate $\hs $ may be consistent. This is similar to our conclusion made above.  The MSE in Cases $2$ and $4$ seems to drop much more rapidly, which reflects a rapid decline in the bias of $\hs$. 
%In fact in Case $4$, the relationship between logarithm of MSE with $\log(s)$ is almost linear, which indicates the bias of $\hs$ decreases at a polynomial rate with $s$.   

\section{Proofs of the Main Results}
\label{sec:proofs}
We start by observing that by Lemma \ref{prp:fund} in Section \ref{sec:technical}
\begin{equation}
\hs=I_1+I_2+I_3+I_4 ,
\label{eq:main2}
\end{equation}
where
\begin{align}
I_1=\frac{1}{\Ns}\sum^{\ns_0}_{j=1}\left(X_{0j}-\mx{0}\right)^2,&~
I_2=\frac{\ns_0}{\Ns}\left(\mx{0}-\hm_0\right)^2\nonumber\\
I_3=\frac{1}{\Ns}\sum^s_{i=1}\sum^{\ns}_{j=1}\left(X_{ij}-\mx{i}\right)^2,&~
I_4=\frac{1}{\Ns}\sum^s_{i=1}\ns\left(\mx{i}-\hm_i\right)^2.\label{eq:defis}
\end{align}

\subsection{Proof of Theorem \ref{lem:fnts}:}
\label{subsec:proof-of-1}
We have two populations and $\Ns=m+m^{\prime}s$. Since $\ns_0=m$ and $\ns_1=m^{\prime}s$, $m/\Ns\rightarrow 0$ and $m^{\prime}s/\Ns\rightarrow 1$ as $s\rightarrow \infty$. 
Further $\mx{0}$ does not depend on $s$ and $X_{1j}$ are i.i.d random variables for $j=1,2,\ldots,m^{\prime}s$, with $E\left(X_{11}\right)=\mu_1$.  So by the \emph{strong law of large numbers} (SLLN) 
\begin{equation*}
\frac{m\mx{0}+\sum^{m^{\prime}s}_{j=1}X_{1j}}{\Ns}\longrightarrow \mu_1\text{~~~~~a.s.}.
\end{equation*}
Now using the fact that $\mx{0}$ does not depend on $s$, it follows that:
\begin{align}
\hm_0&\longrightarrow\min\left(\mx{0},\mu_1\right) \,\,\, \mbox{a.s.} \label{equ:mu0-conv}\\
\hm_1&\longrightarrow\max\left(\min\left(\mx{0},\mu_1\right),\mu_1\right)=\mu_1 \,\,\, \mbox{a.s.}
\label{equ:mu1-conv}
\end{align}

Notice that in this case
$I_1 = \frac{1}{\Ns} \sum^m_{j=1}\left(X_{0j}-\mx{0}\right)^2$ and $\mx{0}$ does not depend on 
$s$ so 
\begin{equation}
\lim_{s \rightarrow \infty} I_1 = 0 \,.
\label{equ:I1}
\end{equation} 
Also $I_2 = \frac{\ns_0}{\Ns}\left(\mx{0}-\hm_0\right)^2$, so by equation (\ref{equ:mu0-conv})
we get 
\begin{equation}
\lim_{s \rightarrow \infty} I_2 = 0 \,.
\label{equ:I2}
\end{equation} 
Further $I_4 = \frac{m's}{\Ns} \left(\mx{1}-\hm_1\right)^2$, so using equation (\ref{equ:mu1-conv})
we get 
\begin{equation}
\lim_{s \rightarrow \infty} I_4 = 0 \,.
\label{equ:I4}
\end{equation} 

Now observe that by standard SLLN 
\begin{equation}
I_3=\frac{1}{\Ns}\sum^{m^{\prime}s}_{j=1}\left(X_{1j}-\mx{1}\right)^2 \longrightarrow \sigma^2 \,\,\, \mbox{a.s.}
\end{equation}
So collecting the terms in \eqref{eq:main2} we get $\hs \longrightarrow \sigma^2$ a.s. proving the strong consistency.

To show the asymptotic normality we consider $\sqrt{\Ns}\left(\hs-\sigma^2\right)$.   
Thus by similar argument as in equations (\ref{equ:I1}) and (\ref{equ:I2}) prove that
\begin{equation}
\lim_{s \rightarrow \infty} \sqrt{\Ns} I_1 = 0 = \lim_{s \rightarrow \infty} \sqrt{\Ns} I_2 \,.
\label{equ:sqrt-I1-I2}
\end{equation}

Further note that
\begin{align}
\mx{1}-\hm_1&=min\left(0,\mx{1}-\hm_0\right)\nonumber\\
&=min\left(0,\mx{1}-min\left(\mx{0},\frac{m\mx{0}+m^{\prime}s\mx{1}}{m+m^{\prime}s}\right)\right)\nonumber\\
&=min\left(0,max\left(\mx{1}-\mx{0},\frac{m}{\Ns}\left(\mx{1}-\mx{0}\right)\right)\right)\nonumber\\
&=\frac{m}{\Ns}\left(\mx{1}-\mx{0}\right)\mathbf{1}_{\left\{\mx{1}\le\mx{0}\right\}}.\label{eq:fourth}
\end{align}
So it follows that
\begin{eqnarray}
\sqrt{\Ns}I_4 & = & \frac{m^{\prime} s}{\sqrt{\Ns}}\left(\mx{1}-\hm_1\right)^2 \nonumber \\
              & = & \frac{m^2 m^{\prime}}{(\Ns)^{5/2}}\left(\mx{1}-\mx{0}\right)^2 
                    \mathbf{1}_{\left\{\mx{1}\le\mx{0}\right\}} \nonumber \\
              & \longrightarrow & 0 \,\,\, \mbox{a.s.} \label{equ:sqrt-I4}
\end{eqnarray}

Now using the standard CLT we get
\begin{equation}
\sqrt{\Ns}\left(I_3-\sigma^2\right)=\sqrt{\Ns}\left(\frac{1}{\Ns}\sum^{m^{\prime}s}_{j=1}\left(X_{1j}-\mx{1}\right)^2-\sigma^2\right) \cd \mbox{N}\left(0, 2\sigma^4 \right)\,.
\label{equ:sqrt-I3}
\end{equation} 

Finally using equations \eqref{equ:sqrt-I1-I2}, \eqref{equ:sqrt-I4} and \eqref{equ:sqrt-I3} we
conclude that
\[
\sqrt{\Ns}\left(\hs-\sigma^2\right)\cd \mbox{N}\left(0,2\sigma^4\right) \,.
\]
\ \ \hfill $\square$

\subsection{\bf Proof of Theorem \ref{thm:general}}
\label{subsec:proof-of-2}
We start by noting that since $\hs$ is the least squared estimator of $\sigma^2$ so
\begin{equation}\label{eq:fund}
\hs \le \frac{1}{\Ns} \sum^s_{i=0}\sum^{\ns_i}_{j=1}\left(X_{ij}-\mu_i\right)^2 \,,
\end{equation}
provided that ${\mathbf \mu}$ satisfies the required tree order restriction. Here we note that
specific constraint such as the tree order restriction is not needed to claim equation 
\eqref{eq:fund}, it will hold for any general constraint under which least square estimator is obtained
as long as ${\mathbf \mu}$ satisfies it.

Now it follows that
\begin{align}
\hs\le&\sigma^2_{(s)}=:\frac{1}{\Ns}\sum^s_{i=0}\sum^{\ns_i}_{j=1}\left(X_{ij}-\mu_i\right)^2\label{eq:sigbound}\\
=&\frac{1}{\Ns}\left\{\sum^s_{i=0}\ns_i\left(\mx{i}-\mu_i\right)^2+\sum^s_{i=0}\sum^{\ns_i}_{j=1}\left(X_{ij}-\mx{i}\right)^2\right\}\nonumber
\end{align}
Since $X_{ij}-\mu_i$ are i.i.d. $\mbox{N}\left(0,\sigma^2\right)$, so by SLLN the first assertion
of the theorem follows.  Furthermore, from the fundamental decomposition,  \eqref{eq:main2} and \eqref{eq:defis} we get
\begin{equation}\label{eq:varbound}
0 \le \left[I_2+I_4\right] \le \frac{1}{\Ns}\sum^s_{i=0}\ns_i\left(\mx{i}-\mu_i\right)^2.
\end{equation}
Note that $\sqrt{\ns_i}\left(\mx{i}-\mu_i\right)\sim \mbox{N}\left(0,\sigma^2\right)$, which implies $\ns_i\left(\mx{i}-\mu_i\right)^2\sim \sigma^2\chi^2_1$.  Thus from the SLLN it follows that $(s+1)^{-1}\sum^s_{i=0}\ns_i\left(\mx{i}-\mu_i\right)^2\rightarrow 2\sigma^2$ as $s\rightarrow\infty$.  

By assumption $(s+1)/\Ns\rightarrow 0$, thus
\[
0 \le \left[I_2+I_4\right]\le\frac{s+1}{\Ns}\frac{1}{s+1}\sum^s_{i=0}\ns_i\left(\mx{i}-\mu_i\right)^2\rightarrow 0 \,\,\, \mbox{a.s.}
\]
Moreover,
\[0\le\hs-\sigma^2_{(s)}=\frac{s+1}{\Ns}\frac{1}{s+1}\sum^s_{i=0}\left(\mx{i}-\mu_i\right)^2-\left[I_2+I_4\right]\rightarrow 0 \,\,\, \mbox{a.s.}
\]
Now using the SLLN as before $\sigma^2_{(s)}\rightarrow\sigma^2$ almost surely. So $\hs$ is strongly consistent.

Now we assume that $s/\sqrt{\Ns}\rightarrow 0$.
To prove the CLT we observe that
\[
\sqrt{\Ns}\left(\hs-\sigma^2\right)=\sqrt{\Ns}\left(I_1+I_3-\sigma^2\right)+\sqrt{\Ns}\left(I_2+I_4\right).
\]
Now from \eqref{eq:varbound}
\[\sqrt{\Ns}\left(I_2+I_4\right)\le\frac{1}{\sqrt{\Ns}}\sum^s_{i=0}\ns_i\left(\mx{i}-\mu_i\right)^2\le\frac{s+1}{\sqrt{\Ns}}\frac{1}{s+1}\sum^s_{i=0}\ns_i\left(\mx{i}-\mu_i\right)^2.
\]
Thus, using similar argument as above we get $s/\sqrt{\Ns}\rightarrow 0 \implies \sqrt{\Ns}\left(I_2+I_4\right)\rightarrow 0$ a.s.

Let us now denote $Y_s=\sum^s_{i=0}\sum^{\ns_i}_{j=1}\left(X_{ij}-\mx{i}\right)^2$.  Note that $\sum^{\ns_i}_{j=1}\left(X_{ij}-\mx{i}\right)^2$ are i.i.d. $\sigma^2\chi^2_{\ns_i-1}$ distributed random variables.  Thus $Y_s\sim\sigma^2\chi^2_{\Ns-s-1}$. Now 
\begin{align}
\sqrt{\Ns}\left(I_1+I_3-\sigma^2\right)&=\sqrt{\Ns}\left(\frac{Y_s}{\Ns}-\sigma^2\right)=\frac{Y_s-E[Y_s]-(s+1)\sigma^2}{\sqrt{\Ns}}\nonumber\\
&=\sqrt{\frac{Var[Y_s]}{\Ns}}\left\{\frac{Y_s-E[Y_s]}{\sqrt{Var[Y_s]}}\right\}-\frac{s+1}{\sqrt{Var[Y_s]}}\sigma^2.\nonumber
\end{align}
Note that from properties of chi-square distribution
\[
\frac{Y_s-E[Y_s]}{\sqrt{Var[Y_s]}}\cd\mbox{N}\left(0,1\right).
\]
Also under our assumption $Var[Y_s]/\Ns=2(\Ns-s-1)\sigma^4/\Ns \rightarrow 2\sigma^4$.  This completes the proof.\hfill$\square$  

\subsection{\bf Proof of Theorem \ref{thm:asf}}
\label{subsec:proof-of-5}
By assumption $\ns_0=\ns=m$ and $\Ns=m(1+s)$ thus from \eqref{eq:defis} it follows that
\[
I_1 + I_3 = \frac{m-1}{m(s+1)}\sum^s_{i=0} \left\{ \sum_{j=1}^m \frac{1}{m-1} \left(X_{ij}-\mx{i}\right)^2 \right\} 
= \frac{m-1}{m} \frac{1}{s+1} \sum_{i=0}^s W_i \,,
\]
where $W_i=\frac{1}{m-1}\sum^m_{j=1}\left(X_{ij}-\mx{i}\right)^2$ for $i \geq 0$.  
Then $\left\{W_i\right\}_{i \geq 0}$ are i.i.d. random variables with mean $\sigma^2$ and variance $2\sigma^4$.
So by standard SLLN we get 
\[
I_1+I_3\rightarrow\frac{m-1}{m}\sigma^2\text{~~~~~as.}
\]

Now consider
\begin{equation}
 I_2 =   \frac{m}{\Ns} \left[\left(\mx{0}-\hm_0\right)^2\right]\nonumber\\
     \le \frac{m}{\Ns} \left(\frac{\rho^{(s)}}{1+\rho^{(s)}}\right)^2
         \left[\left(\lambda(s)\right)^2\right]
     \le \frac{\log s}{1+s} \frac{\left[\left(\lambda(s)\right)^2\right]}{\log s} \,,
\label{eq:t4i2}
\end{equation}
where the first inequality follows from Lemma \ref{lem:meanbound}. Now
\begin{equation}
\lambda(s)=\min_{1\le i\le s}\left(\mx{i}-\mx{0}\right)=\min_{1\le i\le s}\mx{i}-\mx{0}
=\min_{1\le i\le s}\left(\mu_i+\sigma Z_i/\sqrt{m}\right)-\mx{0} \,,
\end{equation}
where $Z_i=\left(\mx{i}-\mu_i\right)/(\sigma/\sqrt{m})$. Note that
$\left\{Z_i\right\}_{i \geq 1}$ are i.i.d $\mbox{N}\left(0,1\right)$ random variables.
From assumption {\bf A1} and {\bf A2} we get
\begin{equation}
\mu_0+\left(\frac{\sigma}{\sqrt{m}}\right)\min_{1\le i\le s}Z_i-\mx{0}\le\lambda(s)\le B+\left(\frac{\sigma}{\sqrt{m}}\right)\min_{1\le i\le s}Z_i-\mx{0}
\label{equ:lambda-bound}
\end{equation}
Observe that $\mx{0}$ does not depend on $s$ and from \citep{gnedenko:hincin:1952}, it follows that $\lambda(s)/\sqrt{2 \log s}\cp \sigma/\sqrt{m}$.  From \eqref{eq:t4i2} it now follows that $I_2\cp 0$. 
 
Finally we consider that
\[
I_4=\frac{1}{1+s} \sum^s_{i=1} \left(\mx{i}-\hm_i\right)^2 \,.
\]
Let $V_{is}=\left(\mx{i}-\hm_i\right)^2$.  From the definition of $\hm_i$, it follows that
\begin{equation} 
\mx{i}-\hm_i\nonumber=\mx{i}-\max\left(\mx{i},\hm_0\right)=\left(\mx{i}-\hm_0\right)\bone_{\left\{\hm_0>\mx{i}\right\}}.
\end{equation} 
From Lemma \ref{lem:meanbound}, $\hm_0\le\mx{0}$, so
\begin{equation} 
\left(\mx{i}-\mx{0}\right)\bone_{\left\{\hm_0>\mx{i}\right\}}\le\left(\mx{i}-\hm_0\right)\bone_{\left\{\hm_0>\mx{i}\right\}}\le0 \,.
\end{equation} 
Thus 
\begin{equation} 
0\le V_{is}=\left(\mx{i}-\hm_0\right)^2 \bone_{\left\{\hm_0>\mx{i}\right\}}
\le \left(\mx{i}-\mx{0}\right)^2 \bone_{\left\{\hm_0>\mx{i}\right\}}
\end{equation} 
Let $\hm_{0(-i)}$ be the estimate of $\mu_0$ obtained after dropping the $i^{\mbox{th}}$ population, that is, using only the data $\left(X_{0k}\right)_{1\le k\le m}$, $\left\{X_{jk}|1\le k\le m\right\}_{1\le j\le s, j\ne i}$.

Recall that $\ns_0=\ns=m$ and notice that
\[
\hm_0\le\hm_{0(-i)}\le\hat{\eta}^{(s)}_{0(-i)}=\min_{1\le j\le s, i\ne j}\left(\frac{\mx{0}+\mx{j}}{2}\right)\text{~~~~~$1\le k\le s$.}  
\]
So it follows that
\begin{align}
V_{is} = & \left(\mx{i}-\mx{0}\right)^2 \bone_{\left\{\hm_0>\mx{i}\right\}} \\
     \le & \left(\mx{i}-\mx{0}\right)^2 \bone_{\left\{\hat{\eta}^{(s)}_{0(-i)}>\mx{i}\right\}}\nonumber\\
      =  & \left(\mx{i}-\mx{0}\right)^2 \mathop{\prod^s_{j=1}}\limits_{j \ne i}
           \bone_{\left\{\frac{1}{2}\left(\mx{0}+\mx{j}\right)>\mx{i}\right\}}\nonumber\\
       = & \left(\mx{i}-\mx{0}\right)^2 \mathop{\prod^s_{j=1}}\limits_{j \ne i}
           \bone_{\left\{\mx{j}>2\mx{i}-\mx{0}\right\}} \label{equ:v1s}
\end{align}

Let us denote $Z^{(s)}_i=2\mx{i}-\mx{0}$ which has $\mbox{N}\left(2\mu_i-\mu_0,5\sigma^2/m\right)$ distribution.  We take expectation on both sides of equation (\ref{equ:v1s}). For the right hand side we first condition on $\mx{i}$ and $\mx{0}$.  After taking expectation over $\mx{i}$ and $\mx{0}$ we get: 
\begin{align}
E\left[V_{is}\right]&\le E\left[\left(\mx{i}-\mx{0}\right)^2\mathop{\prod^s_{j=1}}\limits_{j \ne i}\left\{1-\Phi\left(\frac{\sqrt{m}}{\sigma}\left(Z^{(s)}_i-\mu_j\right)\right)\right\}\right]\nonumber\\
&\le E\left[\left(\mx{i}-\mx{0}\right)^2\left\{1-\Phi\left(\frac{\sqrt{m}}{\sigma}\left(Z^{(s)}_i-B\right)\right)\right\}^{(s-1)}\right].\label{eq:v2s}
\end{align}

\noindent The last inequality holds since $\mu_0\le\mu_j\le B$, for all $1\le j\le s$ by assumption A2.

Now applying the Cauchy-Schwartz inequality on \eqref{eq:v2s} we get
\begin{equation}\label{eq:v3s}
E\left[V_{is}\right]\le\sqrt{E\left[\left(\mx{i}-\mx{0}\right)^4\right]}\sqrt{E\left[\left\{1-\Phi\left(\frac{\sqrt{m}}{\sigma}\left(Z^{(s)}_i-B\right)\right)\right\}^{2(s-1)}\right]}
\end{equation}

Now notice that $E\left[\left(\mx{i}-\mx{0}\right)^4\right]$ does not depend on $s$ and $\mx{i}-\mx{0}$ is stochastically bounded by a $N\left(B-\mu_0,2\sigma^2/m\right)$ random variable.  So $E\left[\left(\mx{i}-\mx{0}\right)^4\right]<C$ for some $C$, for all $i$.  Furthermore, $\sqrt{m}\left(Z^{(s)}_i-B\right)/\sigma\sim \mbox{N}\left(\sqrt{m}(2\mu_i-\mu_0-B)/\sigma,5\right)$ distribution.  Now since $2\mu_i-\mu_0-B\ge\mu_0-B$, it follows that
\begin{equation}\label{eq:v4s}
E\left[V_{is}\right]\le C\sqrt{E\left[\left\{1-\Phi\left(W\right)\right\}^{2(s-1)}\right]},
\end{equation}
where $W\sim \mbox{N}\left(\sqrt{m}(\mu_0-B)/\sigma, 5\right)$ distribution.

Notice that the right-hand side of \eqref{eq:v4s} does not depend on $i$, so
\[E\left[I_4\right]=\frac{1}{s+1}\sum^s_{i=1}E\left[V_{is}\right]\le C\frac{s}{s+1}\sqrt{E\left[\left\{1-\Phi\left(W\right)\right\}^{2(s-1)}\right]}
\]
Now using the \emph{dominated convergence theorem} (DCT) we conclude that $E\left[I_4\right]\rightarrow 0$ as $s \rightarrow \infty$.
This completes the proof.\hfill$\square$

\subsection{\bf Proof of Theorem \ref{thm:hu}}
\label{subsec:proof-of-6}
\noindent $(1)$ This follows by taking expectations on both sides of \eqref{eq:fund} and noting that $\hu$ is unbiased.

\noindent $(2)$ From \eqref{eq:main2} and \eqref{eq:defis} it follows that:
\[\left(\Ns-s-1\right)\hu=\Ns\left(I_1+I_3\right)=\Ns\left(\hs-I_2+I_4\right).
\]
Thus clearly 
\begin{equation}\label{eq:hubound}
\hs-\frac{\Ns-s-1}{\Ns}\hu=I_2+I_4\ge 0.
\end{equation}

\noindent $(3)$ Note that under our assumptions $\hu$ is unbiased and strongly consistent for $\sigma^2$. Under the conditions of Theorem \ref{lem:fnts},  Theorem \ref{thm:general} when $s/\Ns\rightarrow 0$ and Theorem \ref{thm:hmcons} $\hs\cas\sigma^2$ as well.  From these facts the result follows.\hfill$\square$  

%Using \eqref{eq:hubound} we get
%\begin{equation}\label{eq:hshudiff}
%\mid\hs-\hu\mid\le\frac{s+1}{\Ns}\hu+I_2+I_4=\frac{\frac{s+1}{\Ns}}{1-\frac{s+1}{\Ns}}\left(I_1+I_3\right) + I_2+I_4.
%\end{equation}
%Under the conditions of Theorems \ref{lem:fnts}, \ref{thm:general} and \ref{thm:hmcons} it follows that $I_1+I_3\cp\sigma^2$ and $I_2+I_4\cp 0$.  Furthermore, under the conditions of these theorems, $s/\Ns\rightarrow 0$.  
%Thus the LHS of \eqref{eq:hshudiff} converge in probability to $0$.\hfill$\square$
\section{Technical Results on the MLEs of $\mu$ and $\sigma^2$}
\label{sec:technical}
In this section we present some technical results on the constrained MLEs of $\mu$ and $\sigma^2$ which we have used to prove the theorems. Our first result gives an easy but very important decomposition of the MLE $\hs$ of $\sigma^2$ which we refer as fundamental decomposition.  
The proof of this lemma is obvious. So we omit it.

\begin{lemma}\label{prp:fund}
The constrained MLE $\hs$ of $\sigma^2$ admits the following decomposition
\begin{align}
\hs&=\frac{1}{\Ns}\sum^{\ns_0}_{j=1}\left(X_{0j}-\mx{0}\right)^2+\frac{\ns_0}{\Ns}\left(\mx{0}-\hm_0\right)^2\nonumber\\
&+\frac{1}{\Ns}\sum^s_{i=1}\sum^{\ns}_{j=1}\left(X_{ij}-\mx{i}\right)^2+\frac{\ns}{\Ns}\sum^s_{i=1}\left(\mx{i}-\hm_i\right)^2.\label{eq:main}
\end{align} 
\end{lemma}
Note that the first and the third term in \eqref{eq:main} do not involve the order restricted MLEs of the means. 
Form our assumptions, $\Ns$ increases strictly with $s$, so the asymptotic behaviours of these two terms can be determined from classical results such as the SLLN and CLT of i.i.d. random variables 
with finite second moment. 

The next result gives two very useful upper and lower bounds on the MLE $\hm_0$ of $\mu_0$. 
\begin{lemma}\label{lem:meanbound}
Let $\rho^{(s)}:=s\ns/\ns_0$ and $\lambda^{(s)}:=\min_{~1\le i\le s}\left(\mx{i}-\mx{0}\right)$.  Then 
\begin{equation}
\mx{0}+\frac{\rho^{(s)}}{1+\rho^{(s)}}\lambda^{(s)}\bone_{\{\lambda^{(s)}<0\}}\le\hm_0\le\mx{0}.
\label{equ:lb-to-mean-MLE}
\end{equation}
\end{lemma}

\emph{Proof: }
Let $S=\{1,2,\ldots,s\}$. By definition
\[\hm_0=\min_{I\subseteq S}\frac{\ns_0\mx{0}+\ns\sum_{i\in I}\mx{i}}{\ns_0+\ns|I|}=\mx{0}+\min_{I\subseteq S}\frac{\ns\sum_{i\in I}\left(\mx{i}-\mx{0}\right)}{\ns_0+\ns|I|}.
\]
Suppose $\Lambda^{(s)}_i=\mx{i}-\mx{0}$ and $\lambda^{(s)}=\min_{1\le i\le s}\Lambda^{(s)}_i$.  Fix $I\subseteq S$, $I\ne \emptyset$. Then
\[|I|\lambda^{(s)}\le \sum_{i\in I}\Lambda^{(s)}_i\Rightarrow \frac{\ns |I|\lambda^{(s)}}{\ns_0+|I|\ns}\le \frac{\ns \sum_{i\in I}\Lambda^{(s)}_i}{\ns_0+|I|\ns}.\]

Taking minimum on both sides we get:
\begin{equation}\label{eq:uplim}
\min_{I\subseteq S, I\ne \emptyset}\frac{\ns |I|\lambda^{(s)}}{\ns_0+|I|\ns}\le \min_{I\subseteq S, I\ne \emptyset}\frac{\ns \sum_{i\in I}\Lambda^{(s)}_i}{\ns_0+|I|\ns}.
\end{equation}

Note that the function $f(x)=(\ns x c)/(\ns_0+ x \ns)$ is a non-decreasing function if $c>0$ and non-increasing if $c<0$. So in \eqref{eq:uplim}
\[\min_{I\subseteq S, I\ne \emptyset}\frac{\ns |I|\lambda^{(s)}}{\ns_0+|I|\ns}\bone_{\{\lambda^{(s)}<0\}}\ge \frac{s\ns\lambda^{(s)}}{\ns_0+s\ns}\bone_{\{\lambda^{(s)}<0\}}.
\] 

Now using the observation that $\lambda^{(s)}>0$, $\hm_0=\mx{0}$ the inequality follows. \hspace*{1in}\hfill$\square$

We observe that from Lemma \ref{lem:meanbound} it follows
\begin{equation} 
\frac{\rho^{(s)}}{1+\rho^{(s)}} \lambda^{(s)}\bone_{\{\lambda^{(s)}<0\}}\cp 0 
\implies
\mx{0} - \hm_0 \cp 0 \,.
\label{equ:tech-mean-consist}
\end{equation}

From the above lemma following result follows which we present as a stand alone fact.
Note that in \citep[Theorem 2.5]{chaudhuri_perlman_jspi_2007} Chaudhuri and Perlman proved an weaker version using a different technique.
\begin{proposition}
\label{prop:iff}
Suppose $\ns/\log s \longrightarrow \infty$ as $s \rightarrow \infty$ then $\hm_0 \cp \mu_0$ if and only if $\ns_0 \rightarrow \infty$. 
\end{proposition}

\emph{Proof: } Using similar argument which leads to equation \eqref{equ:lambda-bound} we conclude that
if $\ns/log s \longrightarrow \infty$ then 
$\frac{\rho^{(s)}}{1+\rho^{(s)}} \lambda^{(s)}\bone_{\{\lambda^{(s)}<0\}}\cp 0$ 
holds and hence from \eqref{equ:tech-mean-consist} we get
\[
\mx{0} - \hm_0 \cp 0 \,.
\]
The result follows from the fact that $\mx{0} \cp \mu_0$ if and only if $\ns_0 \rightarrow \infty$.
\hfill$\square$
\section*{Acknowledgement} The first author would like to thank the Department of Statistics and Applied Probability of National University of Singapore, for their kind hospitality.  %This work was partially supported by the grant number $R155000111112$ from National University of Singapore.

%\input{paper_Tree-Order-Variance_26OCT2012.bbl}

%\bibliographystyle{unsrtnat}
%\bibliographystyle{gSTA}
%\bibliography{../../../bibs/mypapers,../../../bibs/treeord}
%\input{gSTA.bib}
\end{document}